\documentclass[red,11pt,a4paper]{article}

\usepackage{cite}

\usepackage{amsmath}
\usepackage{amscd}
\usepackage{amssymb}
\usepackage{latexsym}
\usepackage{color}
\usepackage[normalem]{ulem}

\newtheorem{theorem}{Theorem}[section]
\newtheorem{definition}[theorem]{Definition}
\newtheorem{remark}{Remark}
\newcommand{\ep}{\varepsilon}
\newcommand{\R}{\mathbb{R}}

\newcommand{\vr}{\varrho}
\newcommand{\vre}{\varrho_\ep}
\newcommand{\ue}{u_\ep}
\newcommand{\p}{\partial}
\newcommand{\Grad}{\nabla}
\newcommand{\n}{\nabla}
\newcommand{\pt}{\partial_{t}}
\newcommand{\px}{\partial_x}
\newcommand{\lr}[1]{\left( #1 \right)}
\newcommand{\vp}{\varphi}
\newcommand{\va}{\varphi}
\newcommand{\dx}{{\rm d} {x}}
\newcommand{\dt}{{\rm d} t }
\newcommand{\Dt}{\frac{ d}{dt}}
\newcommand{\dxdt}{\dx \,\dt}
\newcommand{\vc}[1]{{\bf #1}}
\newcommand{\eq}[1]{\begin{equation}
\begin{split}
#1
\end{split}
\end{equation}}
\newcommand{\eqh}[1]{\begin{equation*}
\begin{split}
#1
\end{split}
\end{equation*}}
\newcommand{\bo}{| _{\partial\Omega}}
\newcommand{\intO}[1]{\int_{\Omega} #1 \ \dx}
\newcommand{\intTO}[1]{\int_0^T\!\!\!\! \int_{\Omega} #1 \ \dxdt}
\newcommand{\Ov}[1]{\overline{#1}}
\newcommand{\D}{\p_{xx}}

\title{From the highly compressible Navier-Stokes equations to the Porous Medium equation -- rate of convergence}

\author{Boris Haspot 
\thanks{Ceremade UMR CNRS 7534
Universit\'e  Paris  Dauphine,
Place du MarÃÂchal DeLattre De Tassigny
75775 PARIS CEDEX 16 , haspot@ceremade.dauphine.fr}, ~~~~Ewelina Zatorska \thanks{Institute of Applied Mathematics and Mechanics University of Warsaw
ul. Banacha 2, 02-097 Warszawa} \thanks{Institute of Mathematics,
Polish Academy of Sciences,
ul. \'Sniadeckich 8, 
00-656 Warszawa, 
e.zatorska@mimuw.edu.pl}}
\begin{document}
 \maketitle 
\begin{center}
{\bf Abstract}
\end{center}
We consider the one-dimensional Cauchy problem for the Navier-Stokes equations with degenerate viscosity coefficient in highly compressible regime. It corresponds to the compressible Navier-Stokes system with  large Mach number equal to $\ep^{-1/2}$ for $\ep$ going to $0$.  When the initial velocity is related to the gradient of the initial density, {the densities solving the compressible Navier-Stokes equations} --$\rho_\ep$ converge to the unique solution to the porous medium equation \cite{PAM,PAM1}. For viscosity coefficient  $\mu(\rho_\ep)=\rho_\ep^\alpha$ with $\alpha>1$, we obtain a rate of convergence of $\rho_\ep$ in $L^\infty(0,T; H^{-1}(\R))$; for $1<\alpha\leq\frac{3}{2}$ the solution $\rho_\ep$ converges in $L^\infty(0,T;L^2(\R))$. For compactly supported initial data, we prove that most of the mass corresponding to solution $\rho_\ep$ is located in the support of the solution to the porous medium equation. The mass outside this support is small in terms of $\ep$.

\section{Introduction}

The compressible Navier-Stokes equations in the multidimensional case with the Mach number equal to $\ep^{-1/2}$ read:
\begin{equation}
\begin{cases}
\begin{aligned}
&\p_t \vr_\ep+{\rm div}(\vr_\ep u_\ep)=0,\\
&\p_t (\vr_\ep u_\ep)+{\rm div}(\vr_\ep u_\ep\otimes u_\ep)-{\rm div}( 2\mu(\vr_\ep) {\rm D} (u_\ep))
-\nabla(\lambda(\vr_\ep){\rm div}u_{\ep})+\ep \n P(\vr_\ep)=0,\\
&\vre(0,x)=\vr_0(x),\quad \vre u_\ep(0,x)=m_0(x),
\end{aligned}
\end{cases}
\label{main0}
\end{equation}
where $\vr_\ep=\vr_\ep(t,x)$ and $u_\ep=u_\ep(t,x)$ denote  the unknown density and the velocity vector field, respectively, $P(\vr_\ep)=\vre^\gamma$, $\gamma>1$, denotes the pressure, ${\rm D}(u_\ep)=\frac{\Grad u_\ep+\Grad^t u_\ep}{2}$ denotes the symmetric part of the gradient of $u_\ep$, $\mu(\vr_\ep)$ and $\lambda(\vr_\ep)$ denote the two Lam\'e viscosity coefficients satisfying
$$\mu(\vre)>0, \quad 2\mu(\vre)+N\lambda(\vre)\geq0,$$
where $N$ is the space dimension.

\medskip

Our purpose in this paper is to study the asymptotic behaviour of the global weak solutions to the one-dimensional Cauchy problem for \eqref{main0}, when $\ep$ goes to $0$, which corresponds to  the highly compressible limit.  It is the opposite to the low Mach number limit leading to the incompressible Navier-Stokes equations. The {latter} has been intensively studied, especially in the case of constant viscosity coefficients. The first works  in the framework of global weak solutions are due to {P.-L. Lions and N. Masmoudi \cite{LM98} and B. Desjardins and E. Grenier \cite{DG99}} (see also \cite{DGLM99}).  The case of global strong solution with small initial data {was considered by S. Klainerman and A. Majda in \cite{KM81} for well-prepared initial data} and more recently in critical spaces for {scaling of the equations} by R. Danchin {for ill-prepared initial data} in \cite{D02}. Up to our knowledge the highly compressible regime has not been so well studied. The only results we are aware of concern Euler equations, see for example \cite{C}, \cite{GuJu}.
The main difficulty to pass to the limit $\ep\to 0$ in \eqref{main0} with constant viscosity coefficients is due to lack of uniform estimates on $\{\vr_\ep\}_{\ep>0}$ (only the $L^1$ norm is conserved). In this case,  it is  not clear whether the sequence of weak solutions to \eqref{main0} converges to solution of the pressureless system. 
However, in the case of  viscosities satisfying particular algebraic relation 
\begin{equation}
\lambda(\vre)=2\vr\mu'(\vre)-2\mu(\vre),
\label{BD}
\end{equation}
additional compactness information is available. The first author proved in  \cite{PAM,PAM1} that  there exists a limit of the sequence $(\rho_\ep,u_\ep)$,  called {\it a quasi-solution}, satisfying the pressureless system when $N\geq 2$. Moreover, the density of this system corresponds to the solution of the fast diffusion, the heat  or the porous medium equation, depending on the choice of the viscosity coefficients. 
Relation (\ref{BD}) introduced by D. Bresch and B. Desjardins  in \cite{BDa,BD} provides a new entropy structure which ensures estimate on the gradient of the density.  In particular, it gives enough compactness information to pass to the limit when $\ep$ goes to $0$.\\
Stability of  global weak solutions to \eqref{main0} with \eqref{BD} for $\ep$ fixed was proved by A. Mellet and A. Vasseur \cite{MV07} using new energy estimate improving the integrability of the velocity.  A particular choice of viscosity coefficients $\mu(\vr)=\vr$, $\lambda(\vr)=0$  satisfying \eqref{BD} leads to the so-called {\it shallow water system}, for which the proof of existence of global weak solutions  has been recently delivered by A. Vasseur and C. Yu \cite{VY}, the extension to the case \eqref{BD} can be found in \cite{VY1}, see also \cite{GiLaVio,Kor}.   Some ideas on construction of solution using \textit{the cold pressure}  (a pressure that is singular at the vacuum) can be found in \cite{BD1,MPZ, EZ2} and in \cite{BrDeZa} using the notion of the $\kappa$-entropy solutions.

In the following paper, we extend \cite{PAM,PAM1}  by proving rate of convergence of $\vr_\ep$ to a solution of the porous medium equation in a suitable functional framework. We restrict to the one-dimensional Euclidean space $\R$, for which system (\ref{main0}) may be rewritten with a single general viscosity coefficient $\mu(\vr)>0$ in the form
\begin{equation}
\begin{cases}
\begin{aligned}
&\pt \vr_\ep+\px(\vr_\ep u_\ep)=0\\
&\pt\lr{\vr_\ep u_\ep}+\px(\vr_\ep u_{\ep}^2)-\px (\mu(\vr_\ep)\px u_\ep) +\ep \px P(\vr_\ep)=0.
\end{aligned}
\end{cases}
\label{main1}
\end{equation}
Following  \cite{global} (see also \cite{para} and \cite{BrDeZa}), we rewrite  the above system in terms of the {\it effective velocity} 
$$v_\ep=u_\ep+\px\vp(\vre), \quad \mbox{where}  \ \vp'(\vr_\ep)=\frac{\mu(\vr_\ep)}{\vr_\ep^2}$$
which leads to a appearance of a parabolic term in the continuity equation
\begin{equation}
\begin{cases}
\begin{aligned}
&\pt \vr_\ep-\p_{x}\lr{\frac{\mu(\vre)}{\vr_\ep}\p_x\vr_\ep}+\px(\vr_\ep v_\ep)=0,\\
&\pt\lr{\vr_\ep v_\ep}+\px(\vr_\ep u_{\ep} v_\ep) +\ep \px P(\vr_\ep)=0.
\end{aligned}
\end{cases}
\label{main2}
\end{equation}
This new regularising effect was observed in  \cite{para}  in the framework
of strong solution corresponding to initial density far away from the vacuum. The same concept was recently used in \cite{BrGiZa, BrDeZa} to construct regular solutions approximating  global weak solutions to the full low Mach number limit system and to the compressible system.

Our starting observation is that when  $v_\ep(0,\cdot)=0$, the sequence $\{\rho_\ep v_\ep\}_{\ep>0}$ converges to $0$  in a suitable functional framework. It  heuristically implies that the limit $\tilde{\rho}=\lim_{\ep\to0}\rho_\ep$ solves  the fast diffusion, the heat or the porous medium equation
\begin{equation}
\begin{cases}
\begin{aligned}
&\pt \tilde{\rho}-\p_{x}\lr{\frac{\mu(\tilde\rho)}{\tilde\rho}\p_x\tilde\rho}=0,\\
&\tilde\rho(0,x)=\rho_0(x),
\end{aligned}
\end{cases}
\label{main2bb}
\end{equation}
where
\begin{equation}\mu(\tilde\rho)=\tilde\rho^{\alpha}\;\;\mbox{for}\;\;\alpha>0.
\label{loi}
\end{equation}
The rate of convergence of $(\rho_\ep-\tilde\rho )$  can be obtained by employing a duality method  from \cite{Va}, used by J. L. V\'azquez to prove uniqueness of the very weak solutions to the porous medium equation.
  
  In the present paper we restrict to the case {$\alpha>1$} because we are particularly interested in dealing with compactly supported initial data. When 
$0<\alpha<\frac{1}{2}$ the entropy discovered by D. Bresch and B. Desjardins (see  \cite{BD}) allows to bound the density from below. Indeed, in \cite{MV07-1D} A. Mellet and A. Vasseur proved the existence of global strong solution to (\ref{main1})  for initial density far away from the vacuum. The main ingredient of the proof was to use entropy from \cite{BD} to estimate $\p_x(\rho^{\alpha-\frac{1}{2}})$ in $L^\infty(0,T;L^2(\R))$ for all $T>0$. This implies boundedness of  ${\rho}^{-1}$ in $L^\infty(0,T;L^\infty(\R))$ which allows to show sufficient regularity of solution in order to prove the uniqueness. This result has been recently extended by the first author to  the case of the shallow water system ($\alpha=1$)  \cite{global}. The main idea was to prove  boundedness of $v_\ep$ appearing in \eqref{main2} by using the structure of the transport equation for $v_\ep$. It allows to bound ${\rho}^{-1}$  by application of the maximum principle to the continuity equation.

Finally, let us mention several important results concerning existence of global weak and strong solutions to system \eqref{main1} with initial density close to vacuum. Presenting the exhaustive overview of the theory in this field is however beyond of the scope of this article.

The existence of global weak solutions to system (\ref{main1}) was proven by Q. Jiu and Z. Xin in \cite{Jiu}. The main difficulty was to construct approximate regular solutions to (\ref{main1}) which verify  all the entropies used in \cite{MV07}. To construct such an approximation,  Q. Jiu and Z. Xin employed the result on existence of global strong solution on bounded domain $[-M,M]$ with Dirichlet boundary conditions.
The latter result was proven by H-L. Li., J. Li and Z. Xin \cite{LLX} who used energy estimate  to control the Lipschitz norm of the velocity.

There are also several interesting papers on system (\ref{main1}) with  free boundary corresponding to  the interface with the vacuum. The first result in this direction is due to D. Hoff and D. Serre \cite{HS} who proved the existence of global strong solution with discontinuous density at the interface for constant viscosity coefficients. An interesting extension of this result to the Neumann boundary conditions is due to  P. B. Mucha  \cite{Mucha1D} and for the reactive system see \cite{LM}.
S. Jiang, Z. Xin and P. Zhang in \cite{JXZ} obtained similar results in the case of degenerate viscosity coefficients.  The initial density considered in this paper is discontinuous at the interface and  therefore is not compatible with the entropy introduced in \cite{BD}.  We shall discuss this problem with more details later, in Section \ref{Perspectives}.

The free boundary problem with initial density continuously connecting to the vacuum was analysed by T. Yang and H. Zhao in  \cite{YaZa}. They proved the existence of weak solution in a finite time interval in Lagrangian coordinates  by showing that for sufficiently small $T$ the $L^\infty(0,T; L^\infty(\Omega))$ norm of $\rho^{\alpha+1}\p_x u$ is bounded. Unfortunately, this is not sufficient to come back to the  Eulerian description since it is not clear how to solve the ordinary differential equation describing the evolution of the  free boundary. In other words, the change of coordinates may not be well defined, see also \cite{YYZ01,YZ02, VYZ03}.

\medskip

The article is structured in the following way. In Section \ref{Sec:2} we formulate the notion of the weak solution to problem \eqref{main1} for $\ep$ fixed and we state our main results. Next, in Section \ref{Sec:3} we recall the result from \cite{Jiu} on existence of global in time weak solutions and derive the estimates independent of $\ep$. In the second part of this section we recall basic properties of the solution to the one-dimensional porous medium equation taken over from the monograph of J.-L. Vazquez \cite{Va}. In Section \ref{Sec:4} we prove the main results of the paper and in the end, in Section \ref{Perspectives}
 we explain  difficulties arising in the study of the free boundary problems and discuss possible extensions of our results. 

{\bf Notation:} The letter C denotes generic constant, whose value may change from line to line.

\section{Main results}\label{Sec:2}
Below we give the definition of a global weak solution to problem \eqref{main1} which we denote by $(\rho_\ep,\ue)$ or by $(\rho, u)$ when no confusion can arise. The existence of such solution was proven by   Q. Jiu and Z. Xin in \cite{Jiu}. We then formulate our main theorems in which we compare $\rho_\ep$ with $\tilde\rho$ the solution to the porous medium equation. 

Let us start with introducing the hypothesis on the initial data, we assume that

\begin{equation}
\begin{aligned}
&\rho_0\geq 0,\;m_0=0\;\;\mbox{a.e. on}\;\;\{ x\in\R;\,\rho_0(x)=0\},\\
& \rho_0\in L^1(\R)\cap L^\infty(\R), \quad
\px \rho_0^{\alpha-\frac{1}{2}}\in L^2(\R)\\
&\frac{m_0^2}{\rho_0}\in L^1(\R),\quad \frac{|m_0|^{2+\nu}}{\rho_0^{1+\nu}}\in L^1(\R),
\end{aligned}
\label{2.1}
\end{equation}
where $\alpha>\frac{1}{2}$ and $\nu>0$ arbitrary small. 

\begin{definition}\label{Def}
Let $\ep$ be fixed. We say that a pair of functions $(\rho,u)$ is a weak solution to system (\ref{main1}) provided that:
\begin{itemize}
\item The density $\rho\geq 0$ a.e., and the following regularity properties hold
$$
\begin{aligned}
&\rho\in L^\infty(0,T;L^1(\R)\cap L^\gamma(\R))\cap C([0,+\infty), (W^{1,\infty}(\R))'),\\
&\px\rho^{\alpha-\frac{1}{2}}\in L^\infty(0,T;L^2(\R)),\sqrt{\rho}u \in L^\infty(0,T;L^2(\R)),
\end{aligned}
$$
where $(W^{1,\infty}(\R))'$ is the dual space of $W^{1,\infty}(\R)$.
\item For any $t_2\geq t_1\geq 0$ and any $\psi\in C^1( [t_1,t_2]\times \R)$, the continuity equation is satisfied in the following sense:
\begin{equation}
\int_{\R}\rho\psi(t_2)\,\dx-\int_{\R}\rho\psi(t_1)\,\dx=
\int^{t_2}_{t_1}\int_{\R}(\rho\p_t\psi+\rho u\p_x\psi)\,\dxdt.
\label{2.3}
\end{equation}
Moreover, for  $\rho v={\rho u{+}\rho\px\vp(\rho)}$, where $\vp'(\rho)=\frac{\mu(\rho)}{\rho^2}$, the following equality is satisfied
\begin{equation}
\int_{\R}\rho\psi(t_2)\,\dx-\int_{\R}\rho\psi(t_1)\,\dx=
\int^{t_2}_{t_1}\int_{\R}\lr{\rho\p_t\psi+\rho v\p_x\psi-\rho\px\vp(\rho)\px\psi}\,\dxdt.
\label{new2.3}
\end{equation}
\item For any {$\psi\in C^{\infty}_{c}([0,T)\times\R)$} the momentum equation is satisfied in the following sense:
\eq{
&\int_{\R}m_0\psi(0)\, \dx+\int^T_0\int_{\R}\lr{\sqrt{\rho}(\sqrt{\rho}u)\p_t\psi+((\sqrt{\rho}u)^2+\ep\rho^\gamma)\p_x\psi}\,\dxdt\\
&\quad{-\langle\rho^\alpha\p_x u,\p_x\psi\rangle}=0,
\label{2.4}
}
\eq{
&\langle\rho^\alpha\p_x u,\p_x\psi\rangle\\
&\quad=-\int^T_0\int_{\R}\rho^{\alpha-\frac{1}{2}}\sqrt{\rho}u{\p_{xx}}\psi\, \dxdt-\frac{2\alpha}{2\alpha-1}\int^T_0\int_{\R}\p_x(\rho^{\alpha-\frac{1}{2}})\sqrt{\rho} u{\px\psi}\, \dxdt.
\label{2.5}
}
\end{itemize}
\end{definition}
Below we give our first main result on the convergence of $(\rho_\ep, u_\ep)$ to a solution to the associated porous medium equation.

\begin{theorem}
\label{theo1}
Let $\gamma>1$, $\alpha>1$. Moreover, assume that the initial data $(\rho_0,m_0)$ satisfy (\ref{2.1}) and that
\eq{\frac{m_0}{\sqrt{\rho_0}}=-\sqrt{\rho_0}\p_x\va(\rho_0).\label{zero_mass}}
1. System (\ref{main1}) admits a global weak solution $(\rho_\ep,u_\ep)$ in the sense of Definition \ref{Def}. In addition, $\rho_\ep$ converges strongly to $\tilde\rho$ -- the strong solution to the porous medium equation 
\begin{equation}
\begin{cases}
\begin{aligned}
&\pt\tilde\rho-\frac{1}{\alpha}\p_{xx}\tilde\rho^{\alpha}=0,\\
&\tilde\rho(0,x)=\rho_0(x),
\label{porous}
\end{aligned}
\end{cases}
\end{equation}
in the following sense, there exists a constant $C>0$ depending on $\rho_0$ such that
\begin{equation}
\|(\rho-\rho_\ep)(t)\|_{H^{-1}(\R)}\leq C \ep^{\frac{1}{2}}t^{\frac{1}{2}}.
\label{H-1}
\end{equation}
2. For $1<\alpha\leq\frac{3}{2}$ there exists a constant $C>0$ depending on $\rho_0$ such that
\eq{\|(\rho_\ep-\tilde\rho)(t)\|_{L^2(\R)}\leq C\ep^{\frac{1}{4}}t^{\frac{1}{4}}
 ,\label{rate_conv}}
for all $t\geq0$.\\
\end{theorem}

{\begin{remark}
As mentioned above, the existence of weak solutions to system (\ref{main1}) was proven by  Q. Jiu and Z. Xin in \cite{Jiu}. The novelty here is that our definition of weak solution incorporates additional equation \eqref{new2.3}.
\end{remark}}

{\begin{remark}
The initial data from \eqref{zero_mass} can be replaced by a condition which imposes this equality only approximately:
\begin{equation}
\left\|\frac{m_0}{\sqrt{\rho_0}}+\sqrt{\rho_0}\p_x\va(\rho_0)\right\|_{L^2(\Omega)}\leq \eta 
\label{erreur}
\end{equation}
with $\eta>0$. It implies that the inequality (\ref{H-1}) is true up to an error term  coming from (\ref{erreur})
\begin{equation}
\|(\rho-\rho_\ep)(t)\|_{H^{-1}(\R)}\leq C \sqrt{\ep+\eta}\;t^{\frac{1}{2}},
\label{bH-1}
\end{equation}
for a constant $C>0$ depending only on the initial data $\rho_0$. The same applies to inequality (\ref{rate_conv}).
\end{remark}}
\begin{remark}
This result is restricted to the case $1<\alpha\leq\frac{3}{2}$ for technical reasons. This assumption allows to bound $\p_x\tilde{\rho}$ and $\p_x\rho_\ep$ in $L^\infty(0,T; L^2(\R))$; actually we get even more that  $\p_x\rho_\ep^{\alpha-1/2}$, $\p_x\tilde{\rho}^{\alpha-1/2}$ are bounded in $L^\infty(0,T; L^2(\R))$ for any $\alpha>1$ (see also \cite{PAM}).
\end{remark}
The next main result is the following.
\begin{theorem}\label{Th:2}
Let $\gamma>1$, $1<\alpha\leq\frac{3}{2}$. Assume that the initial conditions satisfy \eqref{2.1}  and that  there exist two constants $-\infty<a<b<+\infty$ such that
$${\rm supp}[\rho_0]\subset [a,b].$$
Then, there exist constants $C>0$ and $-\infty<a_1<b_1<\infty$ such that
\begin{equation}
\Omega_t:={\rm supp}[\tilde \rho(t,\cdot)]\subset [a_1- C t^{\frac{1}{\alpha+1}}, b_1+C t^{\frac{1}{\alpha+1}}].
\label{finpropagation}
\end{equation}
Moreover, there exists a constant $C>0$ depending only on $\rho_0$ such that
\begin{equation}
\|\rho_\ep(t)\vc{1}_{\Omega_t^c}\|_{L^1(\R)}\leq C\ep^{\frac{1}{4}}t^{\frac{1}{4} }
(1+t^{\frac{1}{2(\alpha+1)}}),
\label{estmass}
\end{equation}
where $\Omega_t^c$ denotes the complement of $\Omega_t$.\\
\end{theorem}
\begin{remark}
{ The finite speed of propagation of the interface \eqref{finpropagation} is a classical result following from the theory of Porous Medium equation.}
\end{remark}
\begin{remark}
Theorem \ref{theo1} provides a rate of convergence of $(\rho_\ep-\tilde{\rho})$ 
and so it extends the results  of the first author proven in \cite{PAM,PAM1} when $N\geq 2$.
Theorem \ref {Th:2} implies that the mass distributed outside $\Omega_t$--the support of  solution to the porous medium equation, is small of order $\ep^{\frac{1}{4}}$.
The evolution of the support of  $\Omega_t$ is known since the interfaces are described by the Darcy law (see Theorem \ref{Darcy}).  However, $\rho_\ep$ may not be a solution to the free boundary problem. In other words, there is no reason  for ${\rm supp}\rho_\ep(t,\cdot)$ to remain compact for all times $t>0$.  However, the mass that may be spread outside of $\Omega_t$ remains small. In this sense  $\tilde{\rho}(t,\cdot)$ is a good approximation of $\rho_\ep(t,\cdot)$ (for any $t\leq \max(1,{c\ep^{-1/(1+\frac{2}{1+\alpha})}})$ for $c$ small) and we can think about the "quasi finite propagation" of the mass.

\end{remark}
\section{Overview of known results}\label{Sec:3}
This section is devoted to the overview and summary of results that will be used to prove our main theorems.
\subsection{Weak solutions to system \eqref{main1}}
We first recall the result on  global in time existence of weak solutions to system \eqref{main1} with $\ep$ being fixed, whose proof can be found in the paper of Q. Jiu, Z. Xin (\cite{Jiu}, Theorem 2.1.)
\begin{theorem}[Existence of weak solutions]\label{Th:3}
Let $\gamma>1$, $\alpha>1/2$ and let $\ep>0$ be fixed. Assume that the initial conditions satisfy \eqref{2.1}. Then, system (\ref{main1}) possesses a global in time weak solution $(\rho,u)$ in the sense of Definition \ref{Def}. Moreover, this solution satisfies the following inequalities uniformly with respect to $\ep$
\begin{equation}
\begin{aligned}
&\rho\in C((0,T)\times\R),
\end{aligned}
\label{2.6}
\end{equation}
\begin{equation}
\begin{aligned}
&\sup_{t\in[0,T]}\int_{\R}\rho\,\dx+\max_{(t,x)\in[0,T]\times\R}\rho\leq c,
\end{aligned}
\label{2.7}
\end{equation}
\begin{equation}
\begin{aligned}
&\sup_{t\in[0,T]}\int_{\R}\lr{|\sqrt{\rho}u|^2+\p_x( \rho^{\alpha-\frac{1}{2}})+\frac{\ep}{\gamma-1}\rho^{\gamma}}\,\dx\\
&\quad+\int^T_0\int_{\R}\lr{\ep[\p_x(\rho^{\frac{\gamma+\alpha-1}{2}})]^2 +\Lambda^2}\, \dxdt \leq c,
\end{aligned}
\label{2.8}
\end{equation}
where $c$ depends only on the initial data and $\Lambda\in L^2((0,T)\times\R)$ is a function satisfying:
\eq{
&\int^T_0\int_{\R}\Lambda \va\, \dxdt\\
&\quad=-\int^T_0\int_{\R}\rho^{\alpha-\frac{1}{2}}\sqrt{\rho}u\p_x\va\, \dxdt-\frac{2\alpha}{2\alpha-1}\int^T_0\int_{\R}\p_x(\rho^{\alpha-\frac{1}{2}})\sqrt{\rho}u \va\, \dxdt.
\label{2.9}
}
\end{theorem}
{\it Proof.} 
The only element that we have to check is verification of the artificial formulation of the continuity equation \eqref{new2.3}. The rest was already proven by  Q. Jiu, Z. Xin  in \cite{Jiu} (see Theorem 2.1).
Their proof combines the existence of global solutions on bounded domain $[-M,M]$, proven in \cite{LLX} with the diagonal procedure allowing to let $M\to\infty$. 
Therefore, we omit the details of construction of this approximate solution and focus only on justification that the formula \eqref{new2.3} is valid.

Let $\{\rho_\delta, u_\delta\}_{\delta>0}$, $\rho_\delta\geq C(\delta)$ be a sequence smooth approximate solutions to \eqref{main1} {with $\mu(\rho_\delta)$ replaced by $\mu_\delta(\rho_\delta)=\mu(\rho_\delta)+\delta\rho_\delta^\theta$, $\theta\in (0,1/2)$}, constructed in \cite{Jiu} on the bounden interval $\Omega=[-M,M]$ for $M$ large, with the initial conditions
$$\rho_\delta(0)=\rho_{\delta,0}, \quad \rho_\delta(0)u_\delta(0)=m_{\delta,0}$$ 
 and the boundary condition
$$u_\delta\bo=0.$$ 
Assume that $\rho_{\delta,0}$, $m_{\delta,0}$ converge to $\rho_0$, $m_0$ in the following sense
\begin{equation}
\begin{aligned}
&\rho_{\delta,0}\to\rho_0, \quad\text{strongly in } L^1(\Omega), \\
&\px( \rho_{\delta,0}^{\alpha-\frac{1}{2}})\to\px( \rho_0^{\alpha-\frac{1}{2}})\quad\text{strongly in } L^2(\Omega)\\
&\frac{m_{\delta,0}^2}{\rho_{\delta,0}}\to\frac{m_0^2}{\rho_0}\quad\text{strongly in } L^1(\Omega),\\
&\frac{|m_{\delta,0}|^{2+\nu}}{\rho_{\delta,0}^{1+\nu}}\to\frac{|m_0|^{2+\nu}}{\rho_0^{1+\nu}}\in L^1(\R)\quad\text{strongly in } L^1(\Omega).
\end{aligned}
\label{new2.1}
\end{equation}
We show that \eqref{new2.3} is weakly sequentially stable when $\delta\to0$. To this purpose let us recall the a-priori estimates derived in \cite{LLX} following the strategy developed by A. Mellet and A. Vasseur  in \cite{MV07} for the multidimensional case.

First recall that the classical energy balance gives for every $T>0$ the following inequality 
\begin{equation}
\begin{aligned}
&\intO{\lr{\rho_\delta \frac{u_\delta^2}{2}+\frac{\ep}{\gamma-1}\rho_\delta^{\gamma}}(T)}+\intTO{{(\rho_\delta^\alpha+\delta\rho_\delta^\theta)}|\p_x u_\delta|^2}\\
&\hspace{4cm}\leq \intO{\lr{\rho_\delta \frac{u_\delta^2}{2}+\frac{\ep}{\gamma-1}\rho_\delta^{\gamma}}(0)}.
\end{aligned}
\label{energieuni2}
\end{equation}
In addition, the BD entropy gives rise to the estimates
 \begin{equation}
\begin{aligned}
&\intO{\lr{\rho_\delta \frac{(u_\delta+\px{\vp_\delta}(\rho_\delta))^2}{2}+\frac{\ep}{\gamma-1}\rho_\delta^{\gamma}}(T)}\\
&+{\frac{\ep}{\alpha-1}}\intTO{\p_x\rho_\delta^{\gamma}\p_x \rho_\delta^{\alpha-1}}
+{\frac{\ep\delta}{\theta-1}}\intTO{\p_x\rho_\delta^{\gamma}\p_x \rho_\delta^{\theta-1}}\\
&\hspace{4cm}\leq \intO{\lr{\rho_\delta \frac{(u_\delta+\px{\vp_\delta}(\rho_\delta))^2}{2}+\frac{\ep}{\gamma-1}\rho_\delta^{\gamma}}(0)},
\end{aligned}
\label{energieuni}
\end{equation}
{where we denoted $\vp_\delta'(\rho_\delta)=\frac{\mu_\delta(\rho_\delta)}{\rho_\delta^2}=\rho_\delta^{\alpha-2}+\delta\rho_\delta^{\theta-2}$.}
Therefore, we obtain the following bounds
\eq{
&\|\sqrt{\rho_\delta}u_\delta\|_{L^\infty(0,T; L^2(\Omega))}+\|\p_x (\rho_\delta^{\alpha-\frac{1}{2}})\|_{L^\infty(0,T;L^2(\Omega))}
{+\delta\|\p_x (\rho_\delta^{\theta-\frac{1}{2}})\|^2_{L^\infty(0,T;L^2(\Omega))}}\\
&+\ep\|\rho_\delta\|^\gamma_{L^\infty(0,T; L^\gamma(\Omega))}+
\ep\|\p_x(\rho_\delta^{\frac{\gamma+\alpha-1}{2}})\|^2_{L^2(0,T; L^2(\Omega))}
{+\ep\delta\|\p_x(\rho_\delta^{\frac{\gamma+\theta-1}{2}})\|^2_{L^2(0,T; L^2(\Omega))}}\\
&+\|{(\rho_\delta^{\alpha/2}+\sqrt{\delta}\vr_\delta^{\theta/2})}\px u_\delta\|_{L^2(0,T; L^2(\Omega))}\leq C,
\label{BDrho}}
where the constant is independent of $\delta$.

Following \cite{LLX} we obtain 
\begin{equation}
\|\rho_\delta\|_{L^\infty(0,T; L^\infty(\Omega))}\leq C.
\label{technorm}
\end{equation}
This estimate can be then used to derive the improved estimate for the velocity, exactly as in \cite{MV07}.{ Testing the momentum equation by $u_\delta|u_\delta|^\nu$ for $\nu>0$ sufficiently small, and testing the continuity equation by $\frac{|u_\delta|^{2+\nu}}{2+\nu}$ and summing the obtained expressions we obtain
\eq{\label{testMV}
&\Dt\intO{\frac{\rho_\delta|u_\delta|^{2+\nu}}{2+\nu}}+(\nu+1)\intO{\mu_\delta(\rho_\delta)|u_\delta|^\nu|\px u_\delta|^2}
\\
&\quad=-\ep\intO{\px\rho^\gamma_\delta u_\delta|u_\delta|^\nu}.}
Due to the boundary conditions the r.h.s. can be integrated by parts and estimated as follow
\eqh{&-\ep\intO{\px\rho^\gamma_\delta u_\delta|u_\delta|^\nu}
=\ep(\nu+1)\intO{\rho_\delta^\gamma|u_\delta|^\nu\px u_\delta }\\
&\leq
r\frac{\ep(\nu+1)}{2}\intO{\rho_\delta^\alpha|u_\delta|^\nu|\px u_\delta|^2}+ C_r\frac{\ep(\nu+1)}{2}\intO{\rho_\delta^{2\gamma-\alpha}|u_\delta|^\nu}.
}
Note that due to our assumptions $2\gamma-\alpha>0$.
For $r$ sufficiently small, the first term can be absorbed by the l.h.s. of \eqref{testMV}, while to control the second term, we can use the natural bounds for \eqref{BDrho}  and \eqref{technorm} by choosing $\nu$ such that $\frac{\nu}{2}<2\gamma-\alpha$. Summarising, \eqref{testMV} yields}
\eqh{
\sup_{t\in[0,T]}\|\rho_\delta|u_\delta|^{2+\nu}\|_{L^1(\Omega)}\leq C.
}
Having these estimates, we can follow the compactness arguments from \cite{MV07} and from \cite{LLX} to prove the convergence of the approximate solution $(\rho_\delta, u_\delta)$ to some weak solution $(\rho_\ep, u_\ep)$ as $\delta\to0$. Namely
\eq{\label{c_delta}
&\rho_\delta\to\rho_\ep \quad\text{in } C([0,T]\times\Ov{\Omega}),\\
&\px\lr{\rho_\delta^{\alpha-{1\over2}}}\to\px\lr{\rho_\ep^{\alpha-{1\over2}}}\quad\text{weakly in } L^2(0,T; L^2(\Omega)),\\
&\sqrt{\rho_\delta} u_\delta\to\sqrt{\rho_\ep} u_\ep,\quad \rho_\delta^{\alpha} u_\delta\to\rho_\ep^\alpha u_\ep \quad\text{strongly in } L^{2+{\nu\over2}}(0,T; L^{2+{\nu\over2}}(\Omega)),\\
&\rho_\delta^\alpha\px u_\delta\to\Lambda \quad\text{weakly in } L^2(0,T; L^2(\Omega)),
}
where $\Lambda$ satisfies \eqref{2.9}. 

At this point, sequential stability of equations \eqref{2.3} and \eqref{2.4} is verified.  To justify  the limit passage in \eqref{new2.3}, one has to prove the convergence of the term $\rho_\delta v_\delta$. Observe that we have
\eqh{\rho_\delta v_\delta=\rho_\delta u_\delta+\rho_\delta\px\vp_\delta(\rho_\delta)=
\rho_\delta u_\delta+\frac{2}{2\alpha-1}\rho_\delta^{1\over2}\px\lr{\rho_\delta^{\alpha-{1\over2}}}
+{\frac{2\delta}{2\theta-1}\rho_\delta^{1\over2}\px\lr{\rho_\delta^{\theta-{1\over2}}}}.
}
Due to  \eqref{c_delta} the first term converges to $\rho_\ep u_\ep$ strongly in $L^2(0,T; L^2(\Omega))$, while the second term converges to $\frac{2}{2\alpha-1}\rho_\ep^{1\over2}\px\lr{\rho_\ep^{\alpha-{1\over2}}}$ weakly in $L^2(0,T; L^2(\Omega))$; {the last one converges to 0}. In consequence
\eqh{
\rho_\delta v_\delta\to \rho_\ep v_\ep \quad\text{weakly in } L^2(0,T; L^2(\Omega)),
}
and so one can pass to the limit in \eqref{new2.3}.

This is the final argument to prove the sequential stability of weak solutions  to \eqref{main2} on the bounded domain $\Omega=[-M,M]$ and with the Dirichlet boundary condition for $u_\delta$. In order to let $M\to \infty$, we combine the diagonal procedure with the convergence of the initial data \eqref{new2.1} as it was done in \cite{Jiu}. $\Box$

\subsection{The porous medium equation}\label{Sec:Porous}
In order to understand the qualitative properties of the limit solution to \eqref{main1},  we recall several important features of the porous medium equations. The majority of them is taken from the
excellent books of J. L. V\'azquez  \cite{Va}, \cite{Vaz}. The porous medium equation can be written as follows:
\begin{equation}
\begin{cases}
\begin{aligned}
&\p_t\tilde{\rho}-\p_{xx}\tilde{\rho}^\alpha=0,\\
&\tilde{\rho}(0,x)=\tilde{\rho}_0(x),
\end{aligned}
\end{cases}
\label{PM}
\end{equation}
with $\alpha>1$. 

In the sequel we shall set $Q=(0,+\infty)\times\R$. Let us recall the notion of global strong solution for  the porous medium equation (\ref{PM}) 
(see \cite{Va} Chapter 9 for more details).
\begin{definition}
We say that a function $\tilde{\rho}\in C([0,+\infty), L^{1}(\R))$ positive is a strong $L^{1}$ solution to problem (\ref{PM})  if:
\begin{itemize}
\item $\tilde{\rho}^{\alpha}\in L^{1}_{loc}(0,+\infty; L^{1}(\R))$ and $\tilde{\rho}_{t}, \D \rho^{\alpha}\in L^{1}_{loc}((0,+\infty)\times\R)$
\item $\tilde{\rho}_{t}=\mu\p_{xx}\tilde{\rho}^{\alpha}$ in distribution sense.
\item $\tilde{\rho}(t)\rightarrow\tilde{ \rho}_{0}$ as $t\rightarrow0$ in $L^{1}(\R)$.
\end{itemize}
\end{definition}
The existence of global strong solution to \eqref{PM} is guaranteed by the following theorem (see \cite{Va} page 197).
\begin{theorem}
\label{theo2.4}
Let $\alpha>1$ For every non-negative function $\tilde{\rho}_{0}\in L^{1}(\R)\cap L^{\infty}(\R)$  there exists a unique global  strong solution $\tilde{\rho}\geq 0$ of (\ref{PM}). Moreover, $\p_{t}\tilde{\rho}\in L^{p}_{loc}(Q)$ for $1\leq p<\frac{(\alpha+1)}{\alpha}$ and:
$$
\begin{aligned}
&\p_{t}\tilde{\rho}\geq-\frac{\tilde{\rho}}{(\alpha-1)t}\;\;\;\mbox{in}\;\;{\mathcal D}^{'}(Q),\\
&\|\p_{t}\tilde{\rho}(t)\|_{L^{1}(\R)}\leq \frac{2\|\tilde{\rho}_{0}\|_{L^{1}(\R)}}{(\alpha-1)t}.
\end{aligned}
$$
Let $\rho_{1}$ and $\rho_{2}$ be two strong solutions of (\ref{PM}) in $(0,T)\times\R$ then for every $0\leq \tau<t$ 
\begin{equation}
\|\big(\rho_{1}-\rho_{2}\big)_{+}(t)\|_{L^{1}(\R)}\leq \|\big(\rho_{1}-\rho_{2}\big)_{+}(\tau)\|_{L^{1}(\R)}.
\label{contraction}
\end{equation}
If  $\rho_{1}$ and $\rho_{2}$ are two strong solution with initial data $(\rho_{1})_0$ and $(\rho_{2})_0$, such that $(\rho_{1})_0(x)\leq (\rho_{2})_0(x)$ in $\R$, then $\rho_{1}(t,x)\leq\rho_{2}(t,x)$ for all  $(t,x)\in(0,+\infty)\times\R$.
\end{theorem}
Likewise above, there exists a theory of global weak solutions with initial data being bounded measures. It covers the case of very important case of self similar solutions, the so called {\it Barenblatt solutions} $U(t,x,M)$, of the form
\begin{equation}
U(t,x,M)={t^{\frac{1}{\alpha+1}}F(x t^{\frac{-1}{\alpha+1}})},\quad F(\xi)=(C-\kappa\xi^2)_{+}^{\frac{1}{\alpha-1}},
\label{Barren}
\end{equation}
where $\kappa=\frac{\alpha-1}{2\alpha(\alpha+1)}$ and $C>0$ with $C=cM^{\frac{1}{\gamma}}$ with $\gamma=\frac{\alpha+1}{2(\alpha-1)}$ and $c$ depends only on $\alpha>1$.
\\
The Barenblatt solutions verify  (\ref{PM}) for $t>0$ in the sense of distributions and with initial data $\tilde{\rho}_0=M\delta_0$, where $\delta_0$ is the Dirac mass. \\
Analogue solutions exist also in the case of fast diffusion equations corresponding to (\ref{PM}) with $\alpha<1$.

\begin{remark}
\label{propagation}
The comparison principle from Theorem \ref{theo2.4} ensures a finite speed of propagation of solutions to the porous medium equation with compactly supported initial data, see for example  \cite{CGT,Va}. In other words, the solution to the porous medium equation with compactly supported initial data remains compactly supported all along the time.
Indeed, it suffices to compare such solution with the Barenblatt solutions.\\
\end{remark}
Let us now recall the so called $L^{1}-L^{\infty}$ smoothing effect (see \cite{Vaz} page 202).
\begin{theorem}
Let $\tilde{\rho}_0\in L^1(\R)$. For every $t>0$ we have:
$$\tilde{\rho}(t,x)\leq C\|\tilde{\rho}_{0}\|_{L^{1}(0,T;\R^{N})}^{\sigma} t^{-\beta},$$
with $\sigma=\frac{2}{N(\alpha-1)+2}$, $\beta=\frac{1}{(\alpha-1)+2}$ and $C>0$ depends only on $\alpha$ and $N$. The exponents are sharp.\\
When $\tilde{\rho}_0$ belongs also to $L^\infty(\R)$ the maximum principle holds.
\end{theorem}


Let us denote by $\tilde v$ the pressure in the sense of porous medium equation, i.e.
$$\tilde v=\frac{\alpha}{\alpha-1}\tilde{\rho}^{\alpha-1}, \quad \tilde v_0=\frac{\alpha}{\alpha-1}\tilde{\rho}_0^{\alpha-1}.$$

 We recall an  important theorem (see \cite{Va} page 376) which gives a precise behaviour of the interfaces when $\tilde{\rho}_0$ has compact support.
\begin{theorem}
\label{Darcy}
Let $x=s(t)$ denote the right interface of the solution $\tilde{\rho}$ to \eqref{PM} emanating from the compactly supported initial data $\tilde{\rho}_0\in L^1(\R)$. For all $t>0$ there exist the one side limits
\begin{equation}
D^{-}_{x}\tilde v(t,s(t))=\lim_{x\rightarrow s(t)^{-}} \p_x \tilde v(t,x),\;\;D^{+}s(t)=\lim_{h\rightarrow 0^{+}} \frac{1}{h}[s(t+h)-s(t)].
\label{15.60}
\end{equation}
Moreover, the Darcy law holds in the form
\begin{equation*}
D^{+}_t s(t)=-D^{-}_x \tilde v(t,s(t)).
\end{equation*}
The same result is true for the left interface.
\end{theorem}

\section{Proof of main results}\label{Sec:4}
Below we present the proofs of our main results stated in Theorem \ref{theo1} and Theorem \ref{Th:2}. We start by proving the rate of convergence from \eqref{rate_conv}, the convergence of the approximate solutions  $(\rho_\ep,\ue)$ to the solutions of the corresponding porous medium equation is a consequence of this estimate, the details can be found in \cite{PAM}.

\subsection*{Proof of Theorem \ref{theo1}}
First let us note that thanks to Theorem \ref{Th:3}, the continuity equation of system \eqref{main1} after the change of variables 
$v_\ep=u_\ep+\p_x\va(\rho_\ep)$ with $\p_x\va(\rho_\ep)=\frac{\mu(\rho_\ep)}{\rho_\ep^2}\p_x\rho_\ep$ reads
\begin{equation}
\begin{aligned}
&\p_t\rho_\ep-\frac{1}{\alpha}\p_{xx}\rho^\alpha_\ep+\p_x(\rho_\ep v_\ep)=0,\\
\end{aligned}
\label{5}
\end{equation}
and it is satisfied in the sense of distributions.\\
From \eqref{zero_mass} and \eqref{energieuni} it also follows that
\begin{equation}
\sup_{t\in[0,T]}\|\sqrt{\rho_\ep}v_\ep(t)\|_{L^2(\R)}\leq\frac{\ep^{\frac{1}{2}}}{(\gamma-1)^{\frac{1}{2}}}\|\rho_0\|_{L^\gamma(\R)}^{\frac{\gamma}{2}}\leq C\ep^{1\over2}.
\label{enerate}
\end{equation}
Next, we consider $\tilde \rho$ -- the solution to the corresponding porous medium equation with  the same initial data $\rho_0$ 
\eq{\label{PM0}
\begin{cases}
\begin{aligned}
&\p_t\tilde\rho-\frac{1}{\alpha}\p_{xx}\tilde\rho^\alpha=0,\\
&\tilde\rho(0,\cdot)=\rho_0,
\end{aligned}
\end{cases}
}
whose main properties were recalled in Section \ref{Sec:Porous}.\\
We now set
$$R_\ep=\rho_\ep-\tilde\rho.$$
It follows from \eqref{5} and \eqref{PM0} that $R_\ep$ satisfies the equation
\eq{\label{RE}
\begin{cases}
\begin{aligned}
&\p_t R_\ep-\frac{1}{\alpha}\p_{xx} (\rho_\ep^\alpha-\tilde\rho^{\alpha})+\p_x(\rho_\ep v_\ep)=0,\\
&R_\ep(0,x)=0, \quad x\in\R,
\end{aligned}
\end{cases}
}
at least in the sense of distributions.

Our goal is 
to estimate a relevant norm of $R_\ep$ in terms of $\ep$. To this purpose, we use duality technic in the spirit of J. L. V\'azquez (see \cite{Va} Section 6.2.1). Testing \eqref{RE} by $\psi\in C^\infty_c((0,T]\times\R)$, we obtain
\begin{equation}
\begin{aligned}
&\int^T_0\int_{\R}\big( R_\ep\p_t \psi+\frac{1}{\alpha} (\rho_{\ep}^{\alpha}-\tilde\rho^{\alpha})\D\psi\big) \dxdt\\
&\quad+\int^T_0\int_{\R}(\rho_{\ep} v_{\ep})\p_x\psi\, \dxdt-\int_{\R}(R_\ep\psi)(T)\,\dx=0.
\end{aligned}
\label{crucial}
\end{equation}
Let us now define
\begin{equation*}
a(t,x)=
\begin{cases}
\begin{array}{lll}
\frac{\rho_{\ep}^{\alpha}-\tilde\rho^{\alpha}}{\rho_{\ep}-\tilde\rho}\quad &\text{if}\  &\rho_{\ep}\ne\tilde\rho\\
0&\text{if}\  &\rho_{\ep}=\tilde\rho.
\end{array}
\end{cases}
\end{equation*}
This definition implies in particular that $\rho_{\ep}^{\alpha}-\tilde\rho^{\alpha}=a R_\ep$.
We can hence rewrite (\ref{crucial}) as
\begin{equation}\label{ref}
\begin{aligned}
&\int^T_0\int_{\R}\big( R_\ep(\p_t \psi+\frac{1}{\alpha}a \D\psi\big)\, \dxdt\\
&\quad+\int^T_0\int_{\R}(\rho_{\ep} v_{\ep})\p_x\psi\, \dxdt-\int_{\R}(R_\ep\psi)(T)\,\dx=0.
\end{aligned}
\end{equation}
The next step consists of solving the {dual problem} in the interval $[-M,M]\subset\R$
\begin{equation}
\begin{cases}
\begin{aligned}
&\p_t\psi+\frac{1}{\alpha} a_n\D\psi=0,\quad (t,x)\in[0,T]\ \times\ (-M,M),\\
&\psi(t,-M)=\psi(t,M)=0,\quad t\in[0,T],\\
&\psi(T,x)=\theta(x), \quad x\in(-M,M),
\end{aligned}
\end{cases}
\label{systimp}
\end{equation}
where $\theta\in C^{\infty}_{0}((-M,M))$. Above, $a_n$ is a smooth approximation of $a$ such that $0<\eta\leq a_n\leq K<+\infty$ (it will be precisely defined later on). This assumption guarantees that system (\ref{systimp}) is parabolic  in the reverse time $t'=T-t$ and it admits  unique smooth solution {$\psi$} on $[0,T]$. 

Since $\psi$ is a solution to (\ref{systimp}), we can use \eqref{ref} to obtain
\eqh{\int_{\R}R_\ep(T)\theta\,\dx=\frac{1}{\alpha}\int^T_0\int_{\R} R_\ep(a-a_n)\D\psi\, \dxdt
+\int^T_0\int_{\R}(\rho_{\ep} v_{\ep})\p_x\psi\, \dxdt,}
therefore
\eq{
&\left|\int_{\R} R_\ep(T)\theta\,\dx\right|\\
&\quad\leq \left|\int^T_0\int_{\R} (\rho_{\ep} v_{\ep})\p_x\psi\, \dxdt\right|
+ \int^T_0\int_{\R}|R_\ep|\,|(a-a_n)|\,|\p_{xx}\psi|\,\dxdt.
\label{dualestim}}
We now need to obtain the a priori estimate for  $\p_{xx}\psi$.
To this purpose, we multiply equation (\ref{systimp})  by $\zeta\D\psi$, where $\zeta=\zeta(t)$ is a smooth and positive function such that $\frac{1}{2}\leq \zeta\leq 1$  and $\p_t\zeta\geq c>0$, we get
\begin{equation}
\int_0^T\int_{\R}\p_t\psi\zeta\D\psi\, \dxdt+\frac{1}{\alpha}\int_0^T\int_{\R}\zeta a_n(\D\psi)^2\, \dxdt=0.
\label{energ}
\end{equation}
Integrating the first term by parts and using the fact that $\psi(t,\cdot)\in C^\infty_0([-M,M])$ {with all derivatives belonging to $C_0([-M,M])$} (see the general theory of \cite{LSU}, page 341), we obtain
$$
\begin{aligned}
&\int_0^T\int_{\R} \p_t\psi\zeta\D\psi\, \dxdt=-\int_0^T \int_{\R} \zeta\p_x\psi \p_x\p_t\psi\, \dxdt =- \frac{1}{2}\int_0^T \int_{\R}\p_t (\p_x\psi)^2 \zeta\, \dxdt,\\
&=\frac{1}{2}\int_0^T \int_{\R}(\p_x\psi)^2 \p_t \zeta\, \dxdt -\frac{1}{2}\int_{\R}((\p_x\psi)^2\zeta)(T)\,\dx+\frac{1}{2}\int_{\R}((\p_x\psi)^2\zeta)(0)\,\dx.
\end{aligned}
$$
Hence, it follows that
\begin{equation}
\int_0^T\int_{\R} \p_t\psi\zeta\D\psi\, \dxdt\geq \frac{1}{2}\int_0^T \int_{\R}(\p_x\psi)^2 \p_t \zeta\, \dxdt  -\frac{1}{2}\int_{\R}((\p_x\psi)^2\zeta)(T)\,\dx.
\end{equation}
Plugging it into (\ref{energ}) we obtain
\begin{equation}
\begin{aligned}
 \frac{1}{2}\int_0^T \int_{\R}(\p_x\psi)^2 \p_t \zeta\, \dxdt  -\frac{1}{2}\int_{\R}((\p_x\theta )^2\zeta)(T)\,\dx+\frac{1}{\alpha}\int_0^T\int_{\R}\zeta a_n(\D\psi)^2 \dxdt\leq 0.
\end{aligned}
\label{energ1}
\end{equation}
In particular, the assumptions on $\zeta$ imply that we have
\begin{equation}
\begin{aligned}
&\int^T_0\int_{\R} (\p_x\psi)^2\dxdt+\frac{1}{\alpha}\int^T_0 \int_{\R} a_n(\D\psi)^2 \dxdt\leq c\|\p_x \theta\|_{L^2(\R)}^2.
\end{aligned}
\end{equation}
Coming back to (\ref{dualestim}) and recalling that $\mbox{supp}\,\psi$ and $\mbox{supp}\,\theta$ are  included included in  $[0,T]\times[-M,M]$ and $[-M,M]$, respectively, we obtain
\begin{equation}
\begin{aligned}
&\left|\int_{\R} R_\ep (T)\theta\,\dx\right|\\
&\leq C\|\p_x \theta\|_{L^2(\R)}\lr{ \lr{\int_0^T\int_{-M}^M\frac{[a-a_n|^2}{a_n}|R_\ep|^2\, \dxdt}^{\frac{1}{2}}+\|\rho_\ep v_\ep\|_{L^2(0,T;L^2(\R))}}.
\end{aligned}
\label{def1}
\end{equation}
At this stage following ideas from \cite{Va} we construct the approximation $a_n$ that is sufficiently regular and verifies $\eta\leq a_n\leq K$ for $\eta>0$ small and for $K>\eta$ large enough. This requires two steps of approximation. The first  consists of taking $a_{K,\eta}=\min(K,\max(\eta,a))$ with $0<\eta<K$, in other words
\begin{equation*}
a_{K,\eta}(t,x)=\left\{
\begin{array}{llc}
K &\mbox{on} &\{ K<a(t,x)\},\\
a &\mbox{on} &\{ \eta\leq a(t,x)\leq K\},\\
\eta &\mbox{on} &\{ {a(t,x)<\eta}\}.
\end{array}\right.
\end{equation*}
Since  $a_{K,\eta}$ is $L^\infty((0,T)\times\R)$, we know that  $a_{K,\eta}$ is also in $L^p((0,T)\times (-M,M))$ for $1\leq p\leq +\infty$. Therefore, in the next step, we can consider a smooth approximation $a_n=a_{K,\eta}*\psi_n$ (with $\psi_n$ a standard regularising kernel) such that for $n\to\infty$ $a_n\to a_{K,\eta}$ strongly in $L^p((0,T)\times (-M,M))$ for all $1\leq p<+\infty$. 

It remains to estimate the r.h.s. of (\ref{def1}). Since $a_{K,\eta}\geq \eta$ we deduce that also $a_n\geq \eta$ for any $n$, therefore

\begin{equation}
\begin{aligned}
&\int^T_0\int^M_{-M}\frac{|a-a_n|^2}{a_n}|R_\ep|^2\, \dxdt\leq \frac{1}{\eta}\int^T_0\int^M_{-M}|a-a_n|^2|R_\ep|^2\, \dxdt\\
&\leq \frac{2}{\eta}\lr{\int_0^T\int_{-M}^M |a_{K,\eta}-a_n|^2 |R_\ep|^2\, \dxdt+\int^T_0\int^M_{-M}((a-K)_{+}+\eta)^2|R_\ep|^2\, \dxdt}.
\end{aligned}
\label{estima}
\end{equation}
We call the two integrals on the r.h.s. $I_1$ and $I_2$, respectively. The integrand of $I_2$ is pointwise bounded by
$${a^2|R_\ep|^21_{a>K}+\eta^2 \vc{1}_{a\leq K}|R_\ep|^2= (\tilde\rho^{\alpha}-\rho_\ep^{\alpha})^2 \vc{1}_{a>K}+\eta^2 \vc{1}_{a\leq K}|R_\ep|^2,}$$
thus, since $ (\tilde\rho^\alpha-\rho_\ep^\alpha)^2$ is in $L^1(0,T;L^1(\R))$ (which follows from \eqref{2.7} and \eqref{2.1}), using the dominated convergence theorem we justify that  $\int^T_0\int^M_{-M} (\tilde\rho^{\alpha}-\rho_\ep^{\alpha})^2 1_{a>K}\, \dxdt\to0$  for $K\to+\infty$.
It implies that for $K$ large enough, we have
\begin{equation}
\begin{aligned}
{\int^T_0\int^M_{-M}\lr{(a-K)_{+}+\eta}^2|R_\ep|^2\, \dxdt}&\leq\eta^2+ \eta^2\int^T_0\int^M_{-M}|R_\ep|^2\\
&\leq \eta^2(1+2TM C),
\end{aligned}
\label{estim1}
\end{equation}
{where we used the fact that $R_\ep$ is uniformly bounded  in $L^\infty((0,T)\times \R)$ for $0<T<+\infty$.}\\
To estimate $I_1$, it suffices to observe that $R_\ep^2$ is actually uniformly bounded in $\ep$ in $L^{p'}((0,T)\times \R)$ for any $1\leq p'\leq +\infty$ due to \eqref{2.7} and due to $L^\infty$ bound on the solution to the porous medium equation. 
We now use the fact that $a_n\to a_{K,\eta}$ when $n\to+\infty$ strongly in $L^p((0,T)\times (-M,M))$ to obtain that $\int_0^T\int_{-M}^M|a-a_n|^2|R_\ep|^2 \dxdt\rightarrow_{n\rightarrow+\infty} 0$. Using this, (\ref{estim1}) and (\ref{estima}) we deduce that for $n=n(\eta,K,M,T)$ large enough we have:
\begin{equation}
\begin{aligned}
&{\int^T_0\int^M_{-M}\frac{|a-a_n|^2}{a_n}|R_\ep|^2\,\dxdt\leq  2\eta(2+ MTC).}
\end{aligned}
\label{estimab}
\end{equation}
Coming back to (\ref{def1}) and using the H\"older inequality we obtain {that there exists $C>0$ such that for every $\eta>0$}:
\begin{equation}
\begin{aligned}
&\left|\int_{\R} R_\ep (T)\theta\, \dx\right|\\
&\leq C\|\p_x \theta\|_{L^2(\R)}\big( \eta (1+MT) +T^{1\over2}\|\sqrt{\rho_\ep}\|_{L^{\infty}(0,T;L^\infty(\R))}\|\sqrt{\rho_\ep} v_\ep\|_{L^\infty(0,T; L^2(\R))}\big).
\end{aligned}
\label{def1a}
\end{equation}
From  (\ref{technorm})  and (\ref{enerate}) we thus deduce that
\begin{equation}
\begin{aligned}
\left|\int_{\R} R_\ep (T)\theta\,\dx \right|& \leq C\|\p_x \theta\|_{L^2(\R)}\frac{\|\rho_0\|_{L^\gamma}^{\frac{\gamma}{2}}}{{\sqrt{\gamma-1}}}\ep^{\frac{1}{2}}T^{1\over2}\|\sqrt{\rho_\ep}\|_{L^{\infty}(0,T;L^\infty(\R))}\\
& \leq {C}\|\p_x \theta\|_{L^2(\R)}\ep^{\frac{1}{2}}T^{1\over2}.
\end{aligned}
\label{def1ac}
\end{equation}
To conclude, let us observe that the above inequality holds for  any $\theta$ from $C^\infty_0((-M,M))$ with arbitrary $M>0$ and that the constant $C$ is independent of $M$. Thus, letting $M$ to $+\infty$ we obtain 
\begin{equation}
\|R_\ep(T)\|_{H^{-1}(\R)} \leq C\ep^{\frac{1}{2}}T^{1\over2},
\label{fincrucial}
\end{equation}
for any $T>0$.
\\
\\
Concerning the second part of Theorem \ref{theo1}, for $1<\alpha\leq\frac{3}{2}$,
since $\p_x\rho_\ep^{\alpha-\frac{1}{2}}$ is uniformly bounded in $L^\infty(0,T;L^2(\R))$ and $\rho_\ep$ is uniformly bounded in $L^\infty(0,T; L^\infty(\Omega))$ we deduce that $\p_x\rho_\ep$ is uniformly bounded in $L^\infty(0,T;L^2(\R))$.  Next, by \eqref{fincrucial} and uniqueness of limits also $\p_x\tilde\rho$ is bounded in $L^\infty(0,T;L^2(\R))$ (the same result can be obtained in dimension $N\geq 2$; see \cite{PAM}).
Thus, by interpolation, the uniform boundedness of $\rho_\ep$, $\tilde \rho$ in $L^\infty(0,T;L^2(\R))$ and (\ref{fincrucial}), imply that
\begin{equation}
\|R_\ep(T)\|_{L^2(\R)} \leq C\ep^{\frac{1}{4}}T^{\frac{1}{4}}.  
\label{fincrucial1}
\end{equation}
This finishes the proof of Theorem \ref{theo1}. $\Box$
{Note that for $1<p<\infty$ the rate of convergence in $L^\infty(0,T;L^p(\R))$ can be proven by interpolation.}

\bigskip
\subsection*{Proof of Theorem \ref{Th:2}}
Having proven the above rate of convergence, our next aim will be to provide some more information about the behaviour of $\rho_\ep$ on the support of the solution to the corresponding porous medium equation.

Under assumptions of Theorem \ref{Th:2} and from the properties of the solutions to the one-dimensional porous medium equation presented in Section \ref{Sec:Porous}, we deduce that the ${\rm supp}\, \tilde\rho(t,\cdot)$ remains compact for all time $t>0$ and it satisfies \eqref{finpropagation}. Indeed, since $\mbox{supp}\rho_0$ is included in $[a,b]$, we can consider the Barenblatt solution $U(t_1,M)$ at time $t_1>0$ such that
$$\|\rho_0\|_{L^\infty(\R)}\leq \|U(t_1,M)\|_{L^\infty(\R)}.$$
Last inequality is verified if $M$ is large enough depending on $\|\rho_0\|_{L^\infty}$ and $\frac{1}{t_1}$ (see formula (\ref{Barren})).
Using the maximum principle we deduce that:
$$\|\tilde{\rho}(t)\|_{L^\infty(\R)}\leq \|U(t_1+t,M)\|_{L^\infty(\R)}.$$
In particular $\mbox{supp}\tilde{\rho}(t,\cdot)$ is included in $[a_1-C(t+t_1)^{\frac{1}{\alpha+1}},b_1+C(t+t_1)^{\frac{1}{\alpha+1}}]$ for some constants $-\infty<a_1<b_1<\infty$ and $C>0$.
\\
Now, we may  write the following sequence of equalities and inequalities
\begin{equation*}
\begin{aligned}
\|\rho_\ep(t)\vc{1}_{\Omega_t^c}\|_{L^1(\R)}&\leq\|\rho_0\|_{L^1(\R)}-\|\rho_\ep(t)\vc{1}_{\Omega_t}\|_{L^1(\R)}\quad \text{(from the mass conservation)}\\
&\leq\|\tilde\rho(t)\vc{1}_{\Omega_t}\|_{L^1(\R)}-\|\rho_\ep(t)\vc{1}_{\Omega_t}\|_{L^1(\R)}\\
&\leq \|(\tilde\rho(t)-\rho_\ep(t)) \vc{1}_{\Omega_t}  \|_{L^1(\R)}\\
&\leq  |\Omega_t|^{\frac{1}{2}}\|\tilde\rho(t)-\rho_\ep(t)  \|_{L^2(\R)} \quad\quad\text{(from the H\"older inequality)}\\
&\leq C \ep^{\frac{1}{4}}t^{\frac{1}{4}}(1+t^{\frac{1}{2(\alpha+1)}}) \qquad \text{(from  \eqref{finpropagation} and \eqref{fincrucial1})}.
\end{aligned}
\end{equation*}
This finishes the proof of Theorem \ref{Th:2}.
$\Box$

\section{Where are the interfaces?}\label{Perspectives}

A natural question arising in the analysis of one-dimensional problems with compactly supported initial data and free boundary conditions is propagation of the interface. From the classical theory for porous medium equation (see Section \ref{Sec:Porous})
we know that the interface moves with a finite speed. Moreover, according to
(\ref{estmass})  we can estimate the mass corresponding to $\rho_\ep$ inside the support of  $\tilde{\rho}$. This could suggest that the interfaces of free-boundary Navier-Stokes equations behave similarly to the interfaces of the porous medium equation. However, the issue of existence of global weak solutions to the free-boundary Navier-Stokes equations remains open.

Below, we present an overview of partial results in this topic and explain the main problems. \\

1. A case of  compactly supported initial density with jump discontinuity at the interface was studied in \cite{JXZ}. In this article S. Jiang, Z. Xin, and P. Zhang  proved existence of global strong solution to \eqref{main1} written in Lagrangian coordinates. Extending their solution outside the free-boundary domain by imposing the Rankine-Hugoniot condition on the velocity, allows to control $\p_x u_\ep$ in $L^1(0,T;L^\infty(\R))$. It enables in particular to obtain global strong solution on the whole $\R$. It is also sufficient to define the free boundary which verifies an ordinary differential equation.  In this case one can  pass from the Eulerian to the Lagrangian coordinates and conversely.
 
It would be therefore interesting to compare the density $\rho_\ep$ from \cite{JXZ}  with solution to  the porous medium equation $\tilde \rho$ emanating from the same initial data $\rho_0$. However, since $\rho_0$ is discontinuous, it is not compatible with assumptions $\sqrt{\rho_0}\p_x\va(\rho_0)\in L^2(\R)$ and $\rho_0\in L^1(\R)$. 
As a consequence, we loose the a priori estimate on $\sqrt{\rho_\ep}v_\ep$ \eqref{energieuni} which was necessary to show damping of  $\p_x(\rho_\ep v_\ep)$ in \eqref{main2}.
In fact, the behaviour  of $\rho_\ep$ is completely different than behaviour of $\tilde\vr$.
The discontinuity at the interface of $\rho_\ep$ exists all along the time, while $\tilde\rho$ becomes continuous in arbitrary finite time (see Remark \ref{propagation}).\\

2. The free-boundary problem with compact initial density, continuously connected to the vacuum was studied in \cite{YaZa}. There, T. Yang, Z-a. Yao, and C. Zhu proved  the local in time existence of weak solutions to \eqref{main1} in  Lagrangian coordinates. Although this setting seems to be more adequate for our considerations,  the existence of interface corresponding to this solution is not clear. Indeed, in order to determine the evolution of the interface one needs to come back to the Eulerian coordinates.
Unfortunately, result obtained in \cite{YaZa} does not guarantee boundedness of 
 the gradient of the velocity  in $L^1(0,T;L^\infty(\R))$. This lack of control on $\p_x u_\ep$ appears at the boundary, where the density vanishes.

For the porous medium equation, the behaviour of the interface in the 1D case is well understood since it verifies the Darcy law (see Theorem \ref{Darcy}). It could be also appropriate to replace the usual free boundary condition for compressible Navier-Stokes equations by the Darcy law at the interfaces.
Indeed in this case we could give a sense to $u_\ep(t,s_1(t))$  (where $s_1$ denotes  the right interface) by considering the limit  $\lim_{x\to s(t)^{-}} u_\ep(t,x)$ on the left hand side of the free boundary $s_1$. However with such condition the Lagrangian change of coordinates is a priori not possible.

\section*{Acknowledgments} {We thank the referee for his/her valuable comments which helped us to improve the presentation of the paper.} The second author wishes to express her gratitude to the University of Paris Dauphine and to its Department of Mathematics for the kind hospitality. She was also supported by the  National  Science  Centre,  Poland,  grant 2014/14/M/ST1/00108  and by the fellowship START of the Foundation for Polish Science.


\end{document}